\documentclass[12pt]{amsart}
\usepackage{amsmath,centernot}
\usepackage{amssymb}
\usepackage[hidelinks]{hyperref}
\usepackage{tikz}

\newtheorem{theorem}{Theorem}
\newtheorem{lemma}[theorem]{Lemma}

\newtheorem{corollary}[theorem]{Corollary}

\numberwithin{theorem}{section}
\numberwithin{equation}{section}

\theoremstyle{definition}

\newtheorem{definition}[theorem]{Definition}
\newtheorem{observation}[theorem]{Observation}
\newtheorem{question}[theorem]{Question}

\newtheorem{remark}[theorem]{Remark}

\DeclareMathOperator{\Bias}{Bias}

\newcommand{\supp}{\operatorname{supp}}


\begin{document}

\title[Bohr topology]{Bohr topology and difference sets for some abelian groups}

\author{John T. Griesmer}
\address{Department of Applied Mathematics and Statistics\\ Colorado School of Mines, Golden, Colorado}
\email{jtgriesmer@gmail.com}

\begin{abstract}   For a fixed prime $p$, $\mathbb F_{p}$ denotes the field with $p$ elements, and $\mathbb F_{p}^{\omega}$ denotes the countable direct sum $\bigoplus_{n=1}^{\infty} \mathbb F_{p}$.  Viewing $\mathbb F_{p}^{\omega}$ as a countable abelian group, we construct a set $A\subseteq \mathbb F_{p}^{\omega}$ having positive upper Banach density while the difference set $A-A:=\{a-b:a,b\in A\}$ does not contain a Bohr neighborhood of any $c\in \mathbb F_{p}^{\omega}$.   For $p=2$ we obtain a stronger conclusion: $A-A$ does not contain a set of the form $g+(B-B)$, where $B$ is piecewise syndetic. This construction answers negatively a variant of the following question asked by several authors: if $A\subseteq \mathbb Z$ has positive upper Banach density, must $A-A$ contain a Bohr neighborhood of some $n\in \mathbb Z$?

We also construct sets $S, A\subseteq \mathbb F_{p}^{\omega}$ such that $S$ is dense in the Bohr topology of $\mathbb F_{p}^{\omega}$, $A$ has positive upper Banach density, and $A+S$ is not piecewise Bohr.  For $p=2$ we show that every translate of $S$ is a set of topological recurrence and $A+S$ is not piecewise syndetic.  These constructions answer a variant of a question asked by the author. \end{abstract}

\maketitle

\section{Introduction}

\subsection{Structure of difference sets}  If $G$ is an abelian group and $A, B\subseteq G$, $t\in G$, write $A+B$ for the \emph{sumset} $\{a+b:a\in A, b\in B\}$,  $A-A$ for the \emph{difference set} $\{a-b:a,b \in A\}$, and $t+A$ for the \emph{translate} $\{t+a:a\in A\}$.
A \emph{Bohr set} in $G$ is a set of the form $\{g\in G: |\chi_{j}(g)-1|<\varepsilon, 1\leq j \leq k\}$, where $\chi_{j}:G\to \mathcal S^{1}$ are homomorphisms from $G$ to the circle group $\mathcal S^{1}:=\{z\in \mathbb C:|z|=1\}$.  A \emph{Bohr neighborhood} is set containing a translate of a Bohr set.  The \emph{Bohr topology} on $G$ is the topology whose neighborhoods are the Bohr neighborhoods.

For $A\subseteq G$, let $d^{*}(A)$ denote the upper Banach density of $A$ (see \S \ref{sectionUBD} for definitions).  Several authors (\cite{BergelsonRuzsa, RuzsaBook, GriesmerIsr, HegyvariRuzsa}) have asked the following question.

\begin{question}\label{question:Main}
  Is the following implication true?  If $A\subseteq \mathbb Z$ has $d^{*}(A)>0$, then $A-A$ contains a Bohr neighborhood.
\end{question}
The answer is suspected to be ``no," but the question remains open.  For each prime $p$, we answer the analogous question for $\mathbb F_{p}^{\omega}$, the direct sum of countably many copies of $\mathbb F_{p}$ (the group with $p$ elements). For $p=2$, our construction is simpler and yields a (possibly) stronger result, so we consider that case separately.  A set $B\subseteq G$ is \emph{syndetic} if $G$ is the union of finitely many translates of $B$, and $B$ is \emph{piecewise syndetic} if there is a syndetic set $B'\subseteq G$ such that for all finite $F\subseteq B'$, there exists $t\in G$ such that $t+F\subseteq B$.  This is equivalent to the usual definition\footnote{A set $S\subseteq G$ is \emph{thick} if for every finite $F\subseteq G$, there exists $t\in G$ such that $t+F\subseteq S$. A set $A\subseteq G$ is called piecewise syndetic if there is a thick set $S\subseteq G$ and a syndetic set $A'\subseteq G$ such that $A'\cap S \subseteq A$.  Standard arguments \cite{GrSingle} show this definition is equivalent to the one given above.} of ``piecewise syndetic."

\begin{theorem}\label{thm:Main}
  For all $c<\frac{1}{2}$, there is a set $A\subseteq \mathbb F_{2}^{\omega}$ having $d^{*}(A)>c$ such that $A-A$ does not contain a set of the form $B-B+g$, where $B$ is piecewise syndetic and $g\in \mathbb F_{2}^{\omega}$.  Consequently, $A-A$ does not contain a Bohr neighborhood of any $g\in \mathbb F_{2}^{\omega}$.
\end{theorem}
\begin{remark} We do not know whether piecewise syndeticity of $B\subseteq \mathbb F_{2}^{\omega}$ implies $B-B$ contains a Bohr neighborhood. This is closely related to a question of Katznelson \cite{KatznelsonChromatic}.
\end{remark}

 For all primes $p$ we have the following, which is (possibly) weaker than Theorem \ref{thm:Main} in the case $p=2$.

\begin{theorem}\label{thm:Oddp}
  Let $p$ be prime.  For all $c<\frac{1}{2}-\frac{1}{2p}$ there exists $A\subseteq \mathbb F_{p}^{\omega}$ such that $d^{*}(A)>c$ and $A-A$ does not contain a Bohr neighborhood.
\end{theorem}

  The densities in Theorems \ref{thm:Main} and \ref{thm:Oddp} cannot be improved: if $G$ is any group and $A\subseteq G$ has $d^{*}(A)>\frac{1}{2}$, then $A-A=G$.    For odd $p$, we have the following lemma, which relies on the well-known fact that when $H$ is a finite group and $A\subseteq H$ has $|A|>\frac{1}{2}|H|$, then $A-A=H$.

\begin{lemma}\label{lem:OverCritical}
  If $p$ is odd and $A\subseteq \mathbb F_{p}^{\omega}$ has $d^{*}(A)>\frac{1}{2}-\frac{1}{2p}$, then $A-A=\mathbb F_{p}^{\omega}$.
\end{lemma}

\begin{proof}
  Suppose $A\subseteq \mathbb F_{p}^{\omega}$ has $d^{*}(A)>\frac{1}{2}-\frac{1}{2p}$, and let $x\in \mathbb F_{p}^{\omega}$, $x\neq 0$. We will show that $x\in A-A$. Let $H:=\langle x \rangle$ be the subgroup generated by $x$.  Then $H$ is isomorphic to $\mathbb Z/p\mathbb Z$.  By \cite[Theorem 3.5]{BergelsonGlasscock} there is a $y\in \mathbb F_{p}^{\omega}$ such that $|(A+y)\cap H|\geq d^{*}(A)|H|> (\frac{1}{2}-\frac{1}{2p})p=\frac{p-1}{2}$; fix such a $y$ and let $A_{y}:=(A+y)\cap H$.  No subset of $H$ has cardinality strictly between $\frac{p-1}{2}$ and $\frac{p+1}{2}$, so $|A_{y}| \geq \frac{p+1}{2}>\frac{1}{2}|H|$.  Thus $A_{y}-A_{y}=H$.  Now $(A_{y}-A_{y})\subseteq A-A$, so $H\subseteq A-A$, and in particular, $x\in A-A$.
\end{proof}

For odd $p$ we make no assertion about whether $A-A$ contains a Bohr neighborhood when $A\subseteq \mathbb F_{p}$ has $d^{*}(A)=\frac{1}{2}-\frac{1}{2p}$.

 \subsection{Motivating background} A theorem of F{\o}lner \cite{Folner54} states that when $d^{*}(A)>0$, $A-A$ contains a set of the form $U\setminus Z$, where $U$ is a Bohr neighborhood of $0\in G$ and $d^{*}(Z)=0$.  For $G=\mathbb Z$, K\v{r}\'{i}\v{z}'s construction \cite{Kriz} exhibits a set $A$ having $d^{*}(A)>0$ such that $A-A$ does not contain a Bohr neighborhood of $0$, and in fact $A-A$ does not contain a set of the form $S-S$, where $S\subseteq \mathbb Z$ is piecewise syndetic.   McCutcheon presents a simplification of K\v{r}\'{i}\v{z}'s construction due to Ruzsa in \cite{McCutcheonAlexandria,McCutcheonBook}. Alan Forrest \cite{ForrestThesis} constructed a set $A$ with the same properties in $\mathbb F_{2}^{\omega}$.  Our construction, outlined in \S \ref{sec:outline}, is similar to these previous constructions. See Part II of \cite{RuzsaBook} for more on difference sets and Bohr neighborhoods.

\subsection{Sumsets} We say that $B\subseteq G$ is \emph{piecewise Bohr} if there is a Bohr set $B'\subseteq G$ such that for all finite $F\subseteq B'$, there exists $t\in G$ such that $t+F\subseteq B$.

If $G$ is a countable abelian group and $A, B \subseteq G$ have positive upper Banach density, then $A+B$ is piecewise syndetic.  For $G=\mathbb Z$, this result is due to Renling Jin \cite{JinSP}.  Bergelson, Furstenberg, and Weiss \cite{BFW} strengthened the conclusion from ``piecewise syndetic" to ``piecewise Bohr."  These results have been generalized to other settings -- see  \cite{BeiglbockUltrafilter,BBF,BjorklundFish,Dinasso,DinassoEtAl,GriesmerThesis, GriesmerIsr}.  The proofs and examples in \cite{GriesmerIsr} raised the following question, stated for $G=\mathbb Z$ as \cite[Question 5.1]{GriesmerIsr}.

\begin{question}\label{questionIsr}  Let $G$ be a countable abelian group and $S\subseteq G$. Let $\tilde{S}$  be the closure of $S$ in $bG$, the Bohr compactification\footnote{We will not use the Bohr compactification in this article except in reference to Question \ref{questionIsr} -- see \cite{RudinFourier} for exposition.} of $G$.  Let $m_{bG}$ denote Haar measure in $bG$.  Which, if any, of the following implications are true?
\begin{enumerate}
  \item[1.] If $m_{bG}(\tilde{S}) > 0$ and $d^{*}(A)>0$, then $S+A$ is piecewise syndetic.

  \item[2.] If $m_{bG}(\tilde{S}) > 0$ and $d^{*}(A)>0$, then $S+A$ is piecewise Bohr.

  \item[3.] If $S$ is dense in the Bohr topology of $G$ and $d^{*}(A)>0$, then $d^{*}(S+A)=1$.
\end{enumerate}
\end{question}

The following theorem provides a negative answer to all parts of Question \ref{questionIsr} when $G$ is $\mathbb F_{2}^{\omega}$.  We say that $S\subseteq G$ is \emph{chromatically intersective} if for every partition of $G$ into finitely many sets $A_{1},\dots, A_{r}$, there exists $i\leq r$ such that $(A_{i}-A_{i})\cap S$ contains a nonzero element of $G$.   Equivalently, $S$ is chromatically intersective if $(B-B)\cap S$ contains a nonzero element whenever $B$ is piecewise syndetic -- see \cite{GrSingle} for a proof of this well-known equivalence.

\begin{theorem}\label{thm:PWsynd}
 For all $\varepsilon>0$, there are sets $S, A\subseteq \mathbb F_{2}^{\omega}$ such that $d^{*}(A)> 1-\varepsilon$, every translate of $S$ is chromatically intersective, and $S+A$ is not piecewise syndetic.
\end{theorem}
The condition ``every translate of $S$ is a chromatically intersective" implies $S$ is dense in the Bohr topology, since every Bohr neighborhood contains a translate of a difference set $B-B$, where $B$ is syndetic.

For $G=\mathbb F_{p}^{\omega}$ where $p$ is odd, the following provides a negative answer to parts 2 and 3 of Question \ref{questionIsr}.

\begin{theorem}\label{thm:PWBohr}
Let $p$ be an odd prime.  For all $\varepsilon>0$, there are sets $S,A\subseteq \mathbb F_{p}^{\omega}$ such that $d^{*}(A)>\frac{1}{2}-\frac{1}{2p}-\varepsilon$, $S$ is dense in the Bohr topology of $\mathbb F_{p}^{\omega}$, and $S+A$ is not piecewise Bohr.
\end{theorem}

\subsection{Upper Banach density}\label{sectionUBD}  A \emph{F{\o}lner sequence} for an abelian group $G$ is a sequence of finite subsets $\Phi_{n}\subseteq G$ such that $\lim_{n\to \infty} \frac{|(\Phi_{n}+g)\cap \Phi_{n}|}{|\Phi_{n}|}=1$ for every $g\in G$.

If $\mathbf \Phi = (\Phi_{n})_{n\in \mathbb N}$ is a F{\o}lner sequence for $G$ and $A\subseteq G$,  the \emph{upper density of $A$ with respect to $\mathbf{\Phi}$} is $\bar{d}_{\mathbf{\Phi}}(A):=\limsup_{n\to \infty}\frac{|A\cap \Phi_{n}|}{|\Phi_{n}|}$; we write $d_{\mathbf{\Phi}}(A)$ if the limit exists.  The \emph{upper Banach density} of $A$ is $d^{*}(A):=\sup\{\bar{d}_{\mathbf{\Phi}}(A): \mathbf{\Phi} \text{ is a F{\o}lner sequence}\}$.  Note that for every $A\subseteq G$, there is a F{\o}lner sequence $\mathbf{\Phi}$ such that $d^{*}(A)=d_{\mathbf{\Phi}}(A)$.

If $(H_{n})_{n\in \mathbb N}$ is an increasing sequence of finite subgroups of $G$ such that $G= \bigcup_{n=1}^{\infty} H_{n}$, then $(H_{n})_{n\in \mathbb N}$ is a F{\o}lner sequence for $G$.

While upper Banach density is not finitely additive, it enjoys the following weaker property.

\begin{lemma}\label{lemAdditivity}
  Let $G$ be a countable abelian group, $g\in G$, and $A\subseteq G$. If $A\cap (g+)=\varnothing$, then $d^{*}(A\cup (g+A))=2d^{*}(A)$.
\end{lemma}
We omit the proof, which is a straightforward application of the relevant definitions.

\subsection{Dynamical interpretation}\label{sec:Dynamic} Theorems \ref{thm:Main} and \ref{thm:Oddp} have interpretations in terms of dynamical systems - see \cite{BergelsonMcCutcheonRecurrence} for definition of ``set of topological recurrence" and ``set of measurable recurrence."  Theorems 2.2 and 2.6 of \cite{BergelsonMcCutcheonRecurrence} lead to the following corollary of Theorems \ref{thm:Main} and \ref{thm:Oddp}.

\begin{corollary}\label{cor:dynamical}
There exists $S\subseteq \mathbb F_{2}^{\omega}$ such that every translate of $S$ is a set of topological recurrence, while $S$ is not a set of measurable recurrence.  For every prime $p$, there is a set $S\subseteq \mathbb F_{p}^{\omega}$ such that $S$ is dense in the Bohr topology of $\mathbb F_{p}^{\omega}$ while $S$ is not a set of measurable recurrence.
\end{corollary}
The combinatorial constructions proving Theorems \ref{thm:Main} and \ref{thm:Oddp} lead naturally to $(C,F)$ constructions (as presented in \cite{danilenko}). One could thereby construct explicit examples of measure preserving $\mathbb F_{p}^{\omega}$-systems witnessing Corollary \ref{cor:dynamical}, but we do not pursue this here.

\subsection{Other groups}  The examples of Theorems \ref{thm:Main} and \ref{thm:Oddp} yield similar results for extensions: if $G$ is a countable abelian group and $\rho:G\to \mathbb F_{2}^{\omega}$ is a surjective homomorphism, and $S, A\subseteq \mathbb F_{2}^{\omega}$ is as in Theorems \ref{thm:Main} and \ref{thm:PWsynd}, then setting $S':=\rho^{-1}(S)$ and $A':=\rho^{-1}(A)\subseteq G$, we have that $A'-A'$ does not contain a set of the form $g+B-B$, where $B\subseteq G$ is piecewise syndetic, while every translate of $S'$ is chromatically intersective, and $d^{*}(A')\geq d^{*}(A)$.  Consequently, $\mathbb Z^{\omega}:=\bigoplus_{n=1}^{\infty} \mathbb Z$ also satisfies the conclusion of Theorem \ref{thm:Main}, as $\mathbb F_{2}^{\omega}$ is a quotient of $\mathbb Z^{\omega}$.   Similarly, the conclusions of Theorems \ref{thm:Oddp} and \ref{thm:PWBohr} can be obtained for any group having $\mathbb F_{p}^{\omega}$ as a quotient.  Question \ref{question:Main} and all parts of Question \ref{questionIsr} remain open for all groups not having $\mathbb F_{p}^{\omega}$ as a quotient for some prime $p$.  A particularly interesting case may be $\mathbb Q^{\omega}$, the direct sum of countably many copies of $\mathbb Q$.

We hope that our proofs can be adapted to answer Question \ref{question:Main}; the methods of \cite{GriesmerRRPD} can be viewed as such an adaptation for a related question: ``Is there a set $S\subseteq \mathbb Z$ such that every translate of $S$ is a set of measurable recurrence, while $S$ is not a set of strong recurrence?"

\subsection{Acknowledgements}  Our proofs are simpler than those in previous revisions of this article. While discussing these results with Michael Bj{\"o}rklund and Joel Moreira, Bj{\"o}rklund suggested pursuing $(C,F)$ constructions.  The current presentation is a result of viewing our examples as $(C,F)$ constructions.

\subsection{Outline of the article}

After outlining the main ideas in \S \ref{sec:outline} we prove Theorems \ref{thm:Main} and \ref{thm:PWsynd} in \S\S \ref{sec:Groups} and \ref{section2Proof}.  Much of the proofs of Theorems \ref{thm:Oddp} and \ref{thm:PWBohr} are identical to the preceding proofs; in \S \ref{sec:Oddp} we present complete details where the proofs differ, then summarize the remainder.

\section{Outline of the proof of Theorems \ref{thm:Main} and \ref{thm:PWsynd}}\label{sec:outline}

In \S\ref{sec:Broad} we explain the coarsest description of our approach, in \S\ref{sec:FinitePieces} we mention some of the combinatorial details.

 \subsection{Broad approach}\label{sec:Broad} We prove Theorem \ref{thm:Main}  by constructing sets $A, S\subseteq \mathbb F_{2}^{\omega}$ with the following properties: every translate of $S$ is chromatically intersective, $d^{*}(A)>\frac{1}{2}-\varepsilon$, and $(A+S)\cap A=\varnothing$.   The last condition is equivalent to $(A-A)\cap S=\varnothing$.  Since every translate of $S$ is chromatically intersective, $A-A$ cannot contain a set of the form $g+(B-B)$, where $B$ is piecewise syndetic and $g\in \mathbb F_{2}^{\omega}$.

Theorem  \ref{thm:PWsynd} is proved using the same $A$ and $S$ constructed in the proof of Theorem \ref{thm:Main}: we show that $((A+S)-(A+S))\cap S=\varnothing$, so $A+S$ is not piecewise syndetic.  Choosing $x\in \mathbb F_{2}^{\omega}$ so that $(A+x)\cap A=\varnothing$, we then set $E=(A+x)\cup A$.  Lemma \ref{lemAdditivity} then implies $d^{*}(E)>1-2\varepsilon$, while partition regularity of piecewise syndeticity shows that $E+S$ is not piecewise syndetic.

To construct $A$ and $S$, we specify an increasing sequence of subgroups $G_{n}\subseteq \mathbb F_{2}^{\omega}$ such that $\mathbb F_{2}^{\omega}=\bigcup_{n=1}^{\infty} G_{n}$.  We find $S_{n}\subseteq G_{n}$ such that $S_{n}+g$ is $n$-chromatically intersective for all $g\in G_{n}$, while $S_{n}$  fails to be $(\frac{1}{2}-\varepsilon)$-density intersective for some prescribed $\varepsilon$. Letting $S=\bigcup_{n=1}^{\infty} S_{n}$, we get that every translate of $S$ is chromatically intersective.

We will find $A_{n}\subseteq G_{n}$ such that $|A_{n}\cap G_{n}|>\bigl(\frac{1}{2}-\varepsilon\bigr)|G_{n}|$, while $(A_{n}+S_{m})\cap A_{n}=\varnothing$ for all $n$, $m$.  Furthermore, we will make $A_{n}\subseteq A_{n+1}$ for each $n$.  Setting $A:=\bigcup_{n=1}^{\infty} A_{n}$, we get $d^{*}(A)\geq \frac{1}{2}-\varepsilon$, and $(A+S)\cap A=\varnothing$.

\subsection{Finite pieces}\label{sec:FinitePieces}   We consider finite approximations to the intersectivity properties mentioned in our results.

If $r\in \mathbb N$ and $S\subseteq G$, we say that $S$ is \emph{$r$-chromatically intersective} if $G=A_{1}\cup \cdots \cup A_{r}$ implies $(A_{i}-A_{i})\cap S$ contains a nonzero element of $G$.  If $\delta>0$, we say that $S$ is \emph{$\delta$-density intersective} if $d^{*}(A)>\delta$ implies $(A-A)\cap S$ contains a nonzero element of $G$.

Let $\mathbb F_{2}^{n}$ be the $n$th cartesian power of $\mathbb F_{2}$.  Write an element of $\mathbb F_{2}^{n}$ as $x=(x_{1},\dots,x_{n})$, let $\supp(x):=\{i\leq n: x_{i}=1\}$, and let $|x|=|\supp(x)|$.  Let $S_{n,k}:=\{x\in \mathbb F_{2}^{n}:|x|\geq n-k\}$, and $A_{n,k}:=\{x\in \mathbb F_{2}^{n}: |x|\leq \frac{n}{2}-k\}$.  The following facts yield finite versions of our results.

\begin{enumerate}
  \item[(i)] $|x-y|\leq n-2k$ for all $x,y\in A_{n,k}$.  Hence $(A_{n,k}-A_{n,k})\cap S_{n,k}=\varnothing$.

\item[(ii)]   $|A_{n,k}|\approx \frac{1}{2}|\mathbb F_{2}^{n}|$ when $n$ is much larger than $k$.

\item[(iii)]  If $x\in \mathbb F_{2}^{n}$ then $x+S_{n,k}$ is $k$-chromatically intersective.
 \end{enumerate}

Facts (i) and (ii) are easy to verify and form the core of many constructions in additive combinatorics, going back to K\v{r}\'{i}\v{z}'s's example \cite{Kriz} and Ruzsa's introduction of niveau sets \cite{RuzsaComponents}.  See \cite{WolfPopular} for exposition and an application.  Fact (iii) is a consequence of Lov\'asz's lower bound for chromatic numbers of Kneser graphs \cite{Lovasz}, exploited in \cite{Kriz}.  Ruzsa's subsequent simplification of K\v{r}\'{i}\v{z}'s's construction, presented by McCutcheon in \cite{McCutcheonAlexandria, McCutcheonBook},  and Forrest's version for $\mathbb F_{2}^{\omega}$ in \cite{ForrestThesis} all use these ideas.

 Our construction, roughly, will piece together copies of the $A_{n,k}$ in $\mathbb F_{2}^{\omega}$ to form $A$, and likewise piece together copies of the $S_{n,k}$ to form $S$, so that $(A-A)\cap S=\varnothing$, while every translate of $S$ is chromatically intersective.  One difficulty in this approach is that two sets $R$, $R'$ may fail to be $\delta$-density intersective for $\delta$ near $\frac{1}{2}$, their union $R\cup R'$ might be $\frac{1}{4}$-density intersective. This difficulty is surmounted in \cite{Kriz}, \cite{ForrestThesis}, and \cite{McCutcheonAlexandria,McCutcheonBook}, which construct a set which is chromatically intersective but not density intersective.  The way those constructions maintain non-density intersectivity for $S$ also prevents some translates of $S$ from being chromatically intersective, so we need a slightly different approach.  We describe  our strategy in more detail in \S\ref{sec:Dense} after introducing a presentation of $\mathbb F_{2}^{\omega}$.

\section{The group of \texorpdfstring{$\mathbb F_{2}$}{F2}-valued functions on \texorpdfstring{$[0,1)$}{[0,1)}}\label{sec:Groups}

In this section we state some definitions and conventions needed for the proof of Theorems \ref{thm:Main} and \ref{thm:PWsynd}. We identify a useful presentation of $\mathbb F_{2}^{\omega}$, as it is more difficult to prove Theorem \ref{thm:Main} using the usual presentation; see Remark \ref{rem:Presentation} for elaboration.

\subsection{Presentation of \texorpdfstring{$\mathbb F_{2}^{\omega}$}{F2}.}\label{secPresentation}  Write the elements of $\mathbb F_{2}$ as $0,1$. Let $I:=[0,1)\subseteq \mathbb R$ be the half-open unit interval. For each $n\in \mathbb N\cup \{0\}$, let $G_{n}$ denote the set of functions $g:I\to \mathbb F_{2}$ which are constant on intervals of the form $\bigl[\frac{j}{2^{n}},\frac{j+1}{2^{n}}\bigr)$, $j\in \{0,\dots, 2^{n}-1\}$. Then $G_{n}$ is a group under pointwise addition, isomorphic to $\mathbb F_{2}^{2^{n}}$.  Observing that $G_{n}\subseteq G_{n+1}$ for each $n$, we let $G:=\bigcup_{n\in \mathbb N} G_{n}$.  Then $G$ is a countable abelian group isomorphic\footnote{One can construct the isomorphism by hand, but it suffices to observe that $\mathbb F_{2}^{\omega}$ and $G$ are both countably infinite vector spaces over the finite field $\mathbb F_{2}$, and all such vector spaces are mutually isomorphic.} to $\mathbb F_{2}^{\omega}$.  Our constructions are easier to define in $G$ rather than in the standard presentation of $\mathbb F_{2}^{\omega}$, so from now on we work with $G$.

\begin{observation}\label{obs:UBD} As $G_{n}\subseteq G_{n+1}$ and $G$ is the union of these subgroups, the sequence $(G_{n})_{n\in \mathbb N}$ is a F{\o}lner sequence for $G$.  Consequently,  we have $d^{*}(A) \geq \limsup_{n\to \infty} \frac{|A\cap G_{n}|}{|G_{n}|}$ for every $A\subseteq G$.
\end{observation}

\subsection{Restrictions to subintervals}\label{secCylinders} For $n\in \mathbb N\cup\{0\}$, let $\Omega_{n}:=\{\bigl[\frac{j}{2^{n}}, \frac{j+1}{2^{n}}\bigr): j\in \{0,\dots, 2^{n}-1\}\}$, a collection of subintervals partitioning $I$.

\begin{figure}
  \begin{tikzpicture}

\draw	(0.3,3) node[anchor=north] {$0\in I$}
		(4,3) node[anchor=north] {$\tfrac12$}
        (5,3) node[anchor=north] {$\tfrac58$}
		(6,3) node[anchor=north] {$\tfrac{3}{4}$}
        (8.3,3) node[anchor=north] {$1$}
        (0.3,0) node[anchor=north] {$0\in I$}
        (4,0) node[anchor=north] {$\tfrac12$}
        (8.3,0) node[anchor=north] {$1$}
        (0,4) node[anchor=east]{1}
        (0,3) node[anchor=east]{0}
        (0,1) node[anchor=east]{1}
        (0,0) node[anchor=east]{0}
        (-0.5,3.5) node[anchor=east]{$g$}
        (-0.5,0.5) node[anchor=east]{$\displaystyle g_{\tau}$};
\foreach \x in {1, 3, 5, 6, 12, 13, 15, 16, 18, 19, 20, 21, 23, 24, 25, 28, 29, 30, 32}
{
   \draw[style=thick] ({(\x-1)/4},4) -- ({\x/4},4);
}

\foreach \x in {2, 4, 7, 8, 9, 10, 11, 14, 17, 22, 26, 27, 31}
{
   \draw[style=thick] ({(\x-1)/4},3) -- ({\x/4},3);
}

\foreach \x in {18, 19, 20, 21, 23, 24}
{
   \draw[style=thick] ({\x-17},1) -- ({\x-16},1);
}

\foreach \x in {17, 22}
{
   \draw[style=thick] ({\x-17},0) -- ({\x-16},0);
}

\filldraw[draw=gray, bottom color = gray, top color = gray, opacity = 0.1] (4,3) -- (6,3) -- (6,4) -- (4,4);

\end{tikzpicture}
\caption{A typical $g\in G_{5}$, $\tau=[\tfrac12,\tfrac34)\in \Omega_{2}$, and $g_{\tau}\in G_{3}$.}\label{fig:g}
\end{figure}
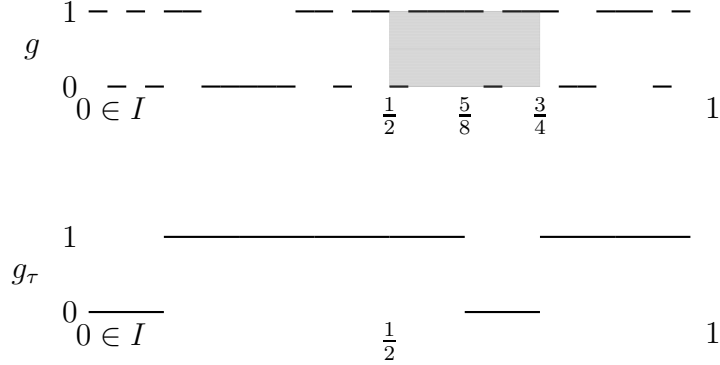

Given $m, n\in \mathbb N$ with $m\leq n$, an interval $\tau=[j2^{-m},(j+1)2^{-m}) \in \Omega_{m}$, and $g\in G_{n}$ we can identify the restriction $g|_{\tau}$ with an element $g_{\tau}\in G_{n-m}$, as shown in Figure \ref{fig:g}.  To be precise: let $\iota_{m,j}:[0,1)\to \tau$, $\iota_{m,j}(t)=(j+t)2^{-m}$, so that $\iota_{m,j}$ is a bijection from $[0,1)$ to $\tau$, mapping intervals $\eta\in \Omega_{n-m}$ of length $2^{-(n-m)}$ onto subintervals of $\tau$ in $\Omega_{n}$ of length $2^{-n}$.  Thus $g|_{\tau}\circ \iota_{m,j}$ is an element of $G_{n-m}$: it is constant on every element of $\Omega_{n-m}$.

For $g\in G_{n}$ and $\tau=[\frac{j}{2^{m}},\frac{j+1}{2^{m}})\in \Omega_{m}$, let $g_{\tau}:=g|_{\tau}\circ \iota_{m,j}$, so that $g_{\tau}\in G_{n-m}$.  For a fixed $\tau\in \Omega_{m}$, and $m\leq n$, the map $g\mapsto g_{\tau}$ is a homomorphism from $G_{n}$ onto $G_{n-m}$.

\subsection{Hamming Balls}\label{sec:HammingBalls}
Let $\mu$ be Lebesgue measure on $[0,1)$.  If $g\in G$, let $\|g\|_{0}:=\mu(g^{-1}(0))$, $\|g\|_{1}:=\mu(g^{-1}(1))$, so that $\|g\|_{0}+\|g\|_{1}=1$ for every $g\in G$.

Let $\mathbf 1 \in G$ denote the constant function where $\mathbf 1(t) = 1 \in \mathbb F_{2}$ for all $t\in I$. Let $\mathbf{0}$ denote the identity element of $G$: $\mathbf 0(t)=0$ for all $t\in I$.

Note that for $g\in G_{n}$, $\|g\|_{0}=2^{-n}\cdot |\{\tau \in \Omega_{n}:g_{\tau}=\mathbf 0\}|$.

For $n, k \in \mathbb N\cup\{0\}$, let
\[
U(n,k):=\{g\in G_{n}:\|g\|_{0}\geq 1-k2^{-n}\}.
\]This is  the \emph{Hamming ball of scale $n$ and radius $k$ around} $\mathbf 0$.  In other words, $U(n,k)$ is the set of $g\in G$ which are constant on intervals of the form $[\frac{j}{2^{n}},\frac{j+1}{2^{n}})$, $0\leq j \leq 2^{n}-1$, and $g_{\tau}=\mathbf{0}$ for at least $2^{n}-k$ such intervals $\tau$. Let $V(n,k) := U(n,k) + \mathbf{1}$, so that
\[V(n,k) = \{g\in G_{n}: \|g\|_{1} \geq 1-k2^{-n}\}.\]
We call $V(n,k)$ the \emph{Hamming ball of scale $n$ and radius $k$ around} $\mathbf 1$.  Note that $\mathbf 1\in V(n,k)$ for all $n,k \in \mathbb N\cup\{0\}$.

\begin{lemma}\label{lem:Vrestriction}
  If $m\leq n$ and $k\in \mathbb N$,
\begin{enumerate}
  \item[(i)] for all $u\in U(n, k)$ and $\tau\in \Omega_{m}$, $u_{\tau}\in U(n-m,k)$.  Similarly if $v\in V(n,k)$, then $v_{\tau}\in V(n-m,k)$.
\item[(ii)] for all $v\in V(n,k)$, all $r,k'\in \mathbb N$,  and at least $n-k$ intervals $\tau\in \Omega_{n}$, $v_{\tau}=\mathbf{1}\in V(r,k')$.   \end{enumerate}

\end{lemma}
\begin{proof}
(i)  If $u\in U(n,k)$, then $u$ is constant on all intervals $\tau \in \Omega_{n}$, and  $u_{\tau}=\mathbf 1$ for at most $k-1$ such intervals.  For an interval $\eta\in \Omega_{m}$, the definition of $u_{\eta}$ then guarantees that $(u_{\eta})_{\theta}=\mathbf 1$  for at most $k-1$ intervals $\theta\in \Omega_{n-m}$, and that $u_{\eta}$ is constant on all such intervals.  Thus $u_{\eta}\in U(n-m,k)$.

Part (ii) follows directly from the definition of $V(n,k)$.
\end{proof}

\begin{remark}
 We call the sets $U(n,k)$ and $V(n,k)$ ``Hamming balls" as we may identify elements of $G_{n}$ with strings of length $2^{n}$ from the alphabet $\mathbb F_{2}$.  With this identification $U(n,k)$ is the set of such strings differing from the constant $0$ string in at most $k$ coordinates.
\end{remark}

\section{Proof of Theorems \ref{thm:Main} and \ref{thm:PWsynd}}\label{section2Proof}

In this section we construct the sets $S, A \subseteq G$ described in Theorem \ref{thm:Main}.   We first show that every translate of a union of the Hamming balls defined in \S\ref{sec:HammingBalls} with increasing scale and radius is chromatically intersective.   We maintain the notation and conventions of \S\ref{sec:outline} and \S\ref{sec:Groups}.

\subsection{Chromatic intersectivity of the \texorpdfstring{$V(n,k)$}{Vnk}}

As in \cite{ForrestThesis, Kriz, McCutcheonAlexandria, McCutcheonBook} we use the following theorem of Lov\'asz \cite{Lovasz}.  See \cite{Matousek} for several proofs and exposition.

\begin{theorem}\label{thmLovasz}
  Let $k, r\in \mathbb N$, and let $E$ be the set of $r$-element subsets of $\{1,\dots, 2r+k\}$.  If $E=\bigcup_{j=1}^{k}E_{j}$, there is a $j\leq k$ and a disjoint pair of elements $e_{1}, e_{2}\in E$ such that $e_{1}, e_{2}\in E_{j}$.
\end{theorem}

For now we fix $n\in \mathbb N$ and work in $\mathbb F_{2}^{n}$, maintaining the notation of \S \ref{sec:FinitePieces}.   Let $k,n\in \mathbb N$, $k\leq n$, and  $H_{n,k}:=\{x\in \mathbb F_{2}^{n}:|x|\leq k\}$.  We will show that every translate of $H_{n,k}$ is $k$-chromatically intersective.  First we prove the following.

\begin{lemma}\label{lem:Poincare}
  If $n \in \mathbb N$,  then $H_{n,2}$ is $n$-chromatically intersective.
\end{lemma}

\begin{proof}
  Write $\mathbb F_{2}^{n}$ as a union $C_{1}\cup \cdots \cup C_{n}$, and consider $E:=\{x\in \mathbb F_{2}^{n}:|x|\leq 1\}$.  Then $|E|=n+1$, so at least one of the $C_{i}$ contains two distinct elements of $E$.  The difference of two elements of $E$ is in $H_{n,2}$, so $(C_{i}-C_{i})\cap H_{n,2}$ contains a nonzero element.
\end{proof}

\begin{lemma}\label{lem:ChromaticBase}
  Fix $n,k\in \mathbb N$, $k\leq n$. Let $S_{n,k}:=\{x\in \mathbb F_{2}^{n}:|x|\geq n-k\}$.  Then $S_{n,k}$ is  $k$-chromatically intersective.
\end{lemma}

\begin{proof}
 For $n=k$ we have $S_{n,k}=\mathbb F_{2}^{n}$, and the conclusion is trivial.  If $n=k-1$ then $S_{n,k}$ contains the nonzero elements of $H_{n,2}$, so the conclusion follows from Lemma \ref{lem:Poincare}.

  Now assume $k\leq n-2$.  Write $n=2r+k$ or $n=2r+k+1$, $r>0$, depending on whether $n-k$ is even or odd.  Write $\mathbb F_{2}^{n}=\bigcup_{i=1}^{k}C_{i}$.  Let $B:=\{x\in \mathbb F_{2}^{n}:|x|=r\}$, and identify $B$ with the collection of $r$-element subsets of $\{1,\dots,n\}$:  $x\leftrightarrow \supp(x)$.  Let $C_{i}':=A\cap C_{i}$, so that the $C_{i}'$ partition $B$.  By Theorem \ref{thmLovasz}, there exists $i\leq k$ and $x, y\in C_{i}$ with $\supp(x)\cap \supp(y)=\varnothing$.  Then $|x-y|=2r\geq n-k$, so $x-y\in S_{n,k}$.
\end{proof}

\begin{lemma}\label{lem:ChromaticTranslate}
Let $n,k\in \mathbb N$, $k\leq n$. For all $y\in \mathbb F_{2}^{n}$, $y+H_{n,k}$ is $k$-chromatically intersective.
\end{lemma}

\begin{proof}
  Let $y\in \mathbb{F}_{2}^{n}$,  let $m=|y|$ and $R:=\supp(y)$.  Consider the group $F_{R}:=\{x\in \mathbb{F}_{2}^{n}: \supp(x)\subseteq R\}$, which is isomorphic to $\mathbb{F}_{2}^{m}$.  We consider three cases.

\noindent \textbf{Case 1: $m\geq k$.} Then $F_{R}\cap (y+H_{n,k})$ contains a copy of $S_{m,k}$.  If $\mathbb F_{2}^{n}$ is partitioned into sets $A_{1}, \dots, A_{k}$, then the sets $A_{i}':=A_{i}\cap F_{R}$ partition $F_{R}$.  Lemma \ref{lem:ChromaticBase} then implies $(A_{i}'-A_{i}')\cap (y+H_{n,k})$ contains a nonzero element for some $i$.

\noindent \textbf{Case 2: $m= k-1$.}  In this case let $R':=R\cup\{t\}$, where $t\notin R$.  Then $(y+H_{n,k})\cap F_{R'}$ contains a copy of $H_{k,2}$, so the conclusion follows from Lemma \ref{lem:Poincare}.

\noindent \textbf{Case 3: $m\leq k-2$.} In this case $y+H_{n,k}$ contains a copy of $H_{n,2}$, and the conclusion again follows from Lemma \ref{lem:Poincare}.
\end{proof}

We now return to the Hamming balls defined in \S\ref{sec:HammingBalls}.

\begin{lemma}\label{lemmaChromatic}
  Let $(n_{i})_{i\in \mathbb N}, (k_{i})_{i\in \mathbb N}$ be sequences of natural numbers with $n_{i},k_{i}\to \infty$, $k_{i}\leq n_{i}$ for each $i$, and let $(g_{i})_{i\in \mathbb N}$ be a sequence of elements of $G$ with $g_{i}\in G_{n_i}$ for each $i$. Then $\bigcup_{i\in \mathbb N} g_{i}+U(n_{i},k_{i})$ is chromatically intersective.  Consequently, every translate of  $S:=\bigcup_{i=1}^{\infty} V(n_{i},k_{i})$ is chromatically intersective.
\end{lemma}

\begin{proof}
  The natural identification of $G_{n}$ with $\mathbb F_{2}^{2^{n}}$ maps $U(n,k)$ to $H_{2^{n},k}$.  Thus Lemma \ref{lem:ChromaticTranslate} implies $g_{i}+U(n_{i},k_{i})$ is $k_{i}$-chromatically intersective.  Fixing $g\in G$, choose $i$ sufficiently large that $g\in G_{n_{i}}$.  Then $g+\mathbf 1\in G_{n_{i}}$, and $g+V(n_{i},k_{i})=g+\mathbf 1+U(n_{i},k_{i})$ is $k_{i}$-chromatically intersective.  Since $k_{i}\to \infty$, we are done.
\end{proof}

\subsection{Some dense subsets of \texorpdfstring{$G$}{G2}}

\subsubsection{Outline}\label{sec:Dense}    Suppose $(n_{i})_{i\in \mathbb N}$, $(k_{i})_{i\in \mathbb N}$ are increasing sequences of integers.  Lemma \ref{lemmaChromatic} shows that every translate of $S:=\bigcup_{i=1}^{\infty} V(n_{i},k_{i})$ is chromatically intersective.  Our goal is to construct a set $A$ having $d^{*}(A)>\frac{1}{2}-\varepsilon$ such that $(A-A)\cap S=\varnothing$; this can be done if $n_{i}\to \infty$ sufficiently rapidly compared to $k_{i}$.  The groups $G_{n}$ form a F{\o}lner sequence, so we find sets $A_{l}\subseteq G_{n_{l}}$ such that $|A_{l}|/|G_{n_{l}}|>\frac{1}{2}-\varepsilon$, while $(A_{l}-A_{l})\cap V(n_{j},k_{j})=\varnothing$ for each $j\leq l$.  The $A_{l}$ are constructed as sumsets with summands given by Definitions \ref{definitionBase} and \ref{def:Citerate}.  In this subsection we indicate the intuition behind those definitions. We first explain how to construct $A_{1}\subseteq G_{n_{1}}$ so that $A_{1}-A_{1}$ is disjoint from $V(n_{1},k_{1})$ and $|A|\approx \frac{1}{2}|G_{n_{1}}|$ (assuming $n_{1}$ is very large compared to $k_{1}$); this is simply a statement of facts (i) and (ii) from \S \ref{sec:FinitePieces} using the presentation of $G$ from \S\ref{sec:Groups}.  Let $A_{1}=\{g\in G_{n_{1}}: \|g\|_{0}\geq \frac{1}{2}+k_{1}2^{-n_{1}}\}$.  Then for $g$, $g'\in A$, we have $\|g-g'\|_{0}\geq 2k_{1}2^{-n_{1}}$, so $g-g'\notin V(n_{1},k_{1})$.  Identifying elements of $G_{n_{1}}$ with subsets of $\{1,\dots, 2^{n_{1}}\}$, it is easy to estimate $|A_{1}|$ using binomial coefficients, and we get that $|A_{1}|\approx \frac{1}{2}|G_{n_{1}}|$ when $n_{1}$ is large compared to $k_{1}$.

We now explain how to construct $A_{2}\subseteq G_{n_{2}}$ such that
\begin{enumerate}
  \item[(i)] $|A_{2}|\approx \frac{1}{2}|G_{n_{2}}|$,
\item[(ii)] $A_{2}-A_{2}$ is disjoint from $V(n_{2},k_{2})$,
\item[(iii)] $A_{2}-A_{2}$ is disjoint from $V(n_{1},k_{1})$. \end{enumerate}
 assuming $n_{2}$ is very large compared to $n_{1}$ and $k_{2}$, and $n_{1}$ is very large compared to $k_{1}$.    After $A_{2}$ is constructed, it is more or less clear how the construction can be iterated to produce the sets $A_{l}$ described above.  However, the arguments to prove properties (i)-(iii) become cumbersome when the obvious definitions are used to perform the iteration.  We take a more streamlined approach in \S\ref{sec:streamline}.

Note that defining $A_{2}$ to be $\{g\in G_{n_{2}}: \|g\|_{0}\geq \frac{1}{2} + k_{2}2^{-n_{2}}\}$ will not satisfy condition (iii) when $n_{2}$ is very large compared to $k_{2}$ and $n_{1}$. Instead, let $m=n_{2}-n_{1}$, and let $A_{2}$ be the set of $g\in G_{n_{2}}$ such that $\|g_{\tau}\|_{0}\geq \frac{1}{2}+k_{2}2^{-m}$ for at least $\frac{1}{2}|\Omega_{n_{1}}|+k_{1}$ intervals $\tau\in \Omega_{n_{1}}$.

\noindent \textit{Proof of (i).} Observe that a typical element of $G_{n_{2}}$ will have, for each $\tau\in \Omega_{n_{1}}$, $\|g_{\tau}\|_{0}\geq \frac{1}{2}+k_{2}2^{-m}$ or $\|g_{\tau}\|_{0}\leq \frac{1}{2}-k_{2}2^{-m}$, and approximately $\frac{1}{2}$ of the elements of $G_{n_{2}}$ will have $\|g_{\tau}\|_{0}\geq \frac{1}{2}+k_{2}2^{-m}$ for at least $\frac{1}{2}|\Omega_{n_{1}}|+k_{1}$  intervals $\tau\in \Omega_{n_{1}}$.

\noindent \textit{Proof of (ii).} For $g$, $g'\in A_{2}$, then there are at least $2k_{1}$ intervals $\tau\in \Omega_{n_{1}}$ such that $\|g_{\tau}\|_{0}, \|g_{\tau}'\|_{0}\geq \frac{1}{2}+k_{2}2^{-m}$, so $g_{\tau}-g_{\tau}'\notin V(m,k_{2})$.   Lemma \ref{lem:Vrestriction} then implies $g-g'\notin V(n_{2},k_{2})$.

\noindent \textit{Proof of (iii).}  The preceding paragraph shows that $g_{\tau}-g_{\tau}'\notin V(n_{2},k_{2})$ for at least $2k_{1}$ intervals $\tau\in \Omega_{n_{1}}$.  In particular $g_{\tau}-g_{\tau}'\neq \mathbf 1$ for at least $2k_{1}$ such intervals. Since $v\in V(n_{1},k_{1})$ has $v_{\tau}=\mathbf 1$ for at least $|\Omega_{n_{1}}|-k_{1}$ such intervals, we conclude that $g-g'\notin V(n_{1},k_{1})$.

To iterate the construction, it is convenient to consider $A_{2}$ as a sumset.  Let $C_{2}:=\{g\in G_{n_{2}}:\|g_{\tau}\|_{0}\geq \frac{1}{2}+k_{2}2^{-m} \text{ for all } \tau \in \Omega_{n_{1}}\}$.  The set $A_{2}$ is not exactly $A_{1}+C_{2}$, but $A_{1}+C_{2}\subset A_{2}$, and $|A_{2}\setminus (A_{1}+C_{2})|$ is very small.

Defining $C_{l}$ as $\{g\in G_{n_{l}}: \|g_{\tau}\|_{0} \geq \frac{1}{2}+k_{l}2^{-(n_{l}-n_{l-1})} \text{ for all } \tau \in \Omega_{n_{l-1}}\}$ (taking $n_{0}=0$, we recover $A_{1}$ as $C_{1}$), we then define $A_{l}$ as a sumset $C_{1}+C_{2}+C_{3}+\cdots + C_{l}$.  The next subsection provides all details.

\subsubsection{Constructing the dense sets}\label{sec:streamline}

\begin{definition}\label{definitionBase} For $n, k \in \mathbb N\cup \{0\}$ let
\[C(0,n;k):= \{g\in G_{n}: \|g\|_{0} \geq \tfrac{1}{2}+k2^{-n}\}.\]
So $C(0,n;k)$ is the set of functions $g:[0,1)\to \mathbb Z/2\mathbb Z$ which are constant on the intervals in $\Omega_{n}$, and the number of such intervals $\tau$ where $g_{\tau}=\mathbf{0}$ is at least $\tfrac{1}{2}|\Omega_{n}|+k$.
\end{definition}

\begin{lemma}\label{lem:uShift}
Let $n, k,k' \in \mathbb N\cup \{0\}$, $k'\leq k$.

\begin{enumerate}
  \item[(i)] $U(n,k')+C(0,n,k) \subseteq C(0,n,k-k')$.
 \item[(ii)] $C(0,n,k)\subseteq C(0,n,k')$.
\end{enumerate}
  \end{lemma}
\begin{proof}
  Part (i) follows from the fact that if $u\in U(n,k)$ and $g\in C(0,n,k)$, then $(u+g)(t)=g(t)$ for all $t\in u^{-1}(0)$.  In other words, $u+g$ and $g$ differ on a set of measure at most $k'2^{-n}$.  Part (ii)  follows immediately from the definition.
\end{proof}

\begin{lemma}\label{lem:BaseDisjoint}
  If $k, k' \geq 0$ and $k+k'>0$, then
\[(C(0,n;k)+\mathbf{1})\cap C(0,n;k')=\varnothing.\]
\end{lemma}
\begin{proof}
  Observe that $(g+\mathbf{1})^{-1}(0)=g^{-1}(1)$.  Thus every $g$ in the intersection satisfies $\|g\|_{0}+\|g\|_{1}\geq \frac{1}{2}+k2^{-n}+ \frac{1}{2}+k'2^{-n}>1$, which is impossible.
\end{proof}

\begin{lemma}\label{lem:BaseCardinality}
  Fix $k\in \mathbb N$.  Then $\lim_{n\to \infty} \frac{|C(0,n;k)|}{|G_{n}|}=\frac{1}{2}$.
\end{lemma}
The conclusion of the lemma can be written as
\begin{equation}\label{eqn:BaseCardinality}
|C(0,n;k)|=\bigl(\tfrac{1}{2}+o(1)\bigr)|G_{n}|,
\end{equation}
where $o(1)$ is a quantity tending to $0$ as $n\to \infty$ (with $k$ fixed).

\begin{proof}
  Identifying an element $g\in C(0,n;k)$ with the collection of intervals $\{\tau \in \Omega_{n}: g_{\tau}= \mathbf{1}\}$, we see that $|C(0,n;k)|$ is equal to the number of subsets of $\{1,\dots,|\Omega_{n}|\}$ having cardinality at most $|\Omega_{n}|-k$.  Thus
\begin{equation*}
    |C(0,n;k)| = \tfrac{1}{2} \Bigl(2^{|\Omega_{n}|}-\sum_{j=-k+1}^{k-1}\binom{|\Omega_{n}|}{\tfrac{1}{2}|\Omega_{n}|-j}\Bigr)\geq \tfrac{1}{2}2^{|\Omega_{n}|}-2k\binom{|\Omega_{n}|}{|\Omega_{n}|/2},
\end{equation*} as $\binom{|\Omega_{n}|}{|\Omega_{n}|/2}\geq \binom{|\Omega_{n}|}{t}$ for all $t$. The standard binomial estimate $\binom{|\Omega_{n}|}{|\Omega_{n}|/2}=o(2^{|\Omega_{n}|})$ yields $C(0,n;k)=(\frac{1}{2}+o(1))2^{|\Omega_{n}|}$.  Since $|G_{n}|=2^{|\Omega_{n}|}$, we are done. \end{proof}

\begin{definition}\label{def:Citerate}
  Fix $m<n\in \mathbb N$ and $k\in \mathbb N\cup\{0\}$.  Let
\begin{align*}
C(m,n;k)&:=\{g\in G_{n}:g_{\tau}\in C(0,n-m;k) \text{ for all } \tau\in \Omega_{m}\}.
\end{align*}
\end{definition}

Elements of $C(m,n;k)$ may be constructed freely by choosing their values on each $\tau$ in $\Omega_{m}$, so
\begin{equation}\label{eqn:Ccount}
|C(m,n;k)|=|C(0,n-m;k)|^{|\Omega_{m}|}.
\end{equation}
Combining (\ref{eqn:Ccount}) with Lemma \ref{lem:BaseCardinality} yields the following estimate.
\begin{lemma}\label{lem:Ccount}
  Fix $m, k\in \mathbb N$. Then
\begin{equation}\label{eqn:Cestimate}
  |C(m,n;k)|=\Big(\frac{1}{|G_{m}|}+o(1)\Bigr)|G_{n}|
\end{equation}
where $o(1)$ is a quantity tending to $0$ as $n\to \infty$ and $m, k$ remain fixed.
\end{lemma}
\begin{proof}
Note that $|G_{n-m}|^{|\Omega_{m}|}=(2^{2^{n-m}})^{(2^{m})}=2^{2^{n}}=|G_{n}|$, while $|G_{m}|=2^{|\Omega_{m}|}$.  Starting with Equation (\ref{eqn:Ccount}), we have
\begin{align*}
  |C(m,n;k)|&=|C(0,n-m;k)|^{|\Omega_{m}|} \\
    &=\Bigl((\tfrac12 + o(1))|G_{n-m}|\Bigr)^{|\Omega_{m}|} && \text{by } (\ref{eqn:BaseCardinality})\\
    &= \Bigl(\frac{1}{2^{|\Omega_{m}|}}+o(1)\Bigr)|G_{n-m}|^{|\Omega_{m}|}\\
    &= \Bigl(\frac{1}{|G_{m}|}+o(1)\Bigr)|G_{n}|.
  \end{align*}
The third line above is obtained by expanding $(\frac{1}{2}+o(1))^{|\Omega_{m}|}$ as $\frac{1}{2^{|\Omega_{m}|}}$ plus a sum of $2^{|\Omega_{m}|}-1$ terms with $o(1)$ as a factor, using the fact that $m$ is fixed.
\end{proof}

\begin{lemma}\label{lem:Uniqueness}
Fix $m\in \mathbb N\cup \{0\}$, $n>m$, and $k\in \mathbb N$.  If $g, g'\in G_{m}$, $h, h'\in C(m,n;k)$ and $g+h=g'+h'$, then $h=h'$ and $g=g'$.
\end{lemma}

\begin{proof}
  Let $g, g' \in G_{m}$, $h,h'\in C(m,n;k)$.  If $g+h=g'+h'$, then $g-g'=h-h'$.  Now $g-g'\in G_{m}$, so $(g-g')_{\tau}\in \{\mathbf 0, \mathbf 1\}$ for all $\tau\in \Omega_{m}$, while the definition of $C(m,n;k)$ implies $h_{\tau},h_{\tau}'\in C(0,n-m;k)$ for all such $\tau$.  Lemma \ref{lem:BaseDisjoint} implies $h_{\tau}-h'_{\tau}\neq \mathbf 1$ for all such $\tau$, so $h_{\tau}-h'_{\tau} = (g-g')_{\tau}= \mathbf 0$ for all $\tau$.  Thus $h=h'$ and $g=g'$.
\end{proof}  With $m,n,k$ as in the lemma, the map $(g,h)\to g+h$ from $G_{m}\times C(m,n;k)$ to $G_{n}$ is one-to-one, and we have the following corollary.
\begin{corollary}\label{cor:Uniqueness}
  If $F\subseteq G_{m}$ then $|F+C(m,n;k)|=|F| |C(m,n;k)|$.
\end{corollary}

\begin{lemma}\label{lem:Cardinality}
Fix $m, k\in \mathbb N$.  If $F\subseteq G_{m}$, then $\lim_{n\to \infty} \frac{|F+C(m,n;k)|}{|G_{n}|}=\frac{|F|}{|G_{m}|}$.
\end{lemma}
The conclusion can be written as $|F+C(m,n;k)|=\bigl(\frac{|F|}{|G_{m}|}+o(1)\bigr)|G_{n}|$, where $o(1)\to 0$ as $n\to \infty$ (and $m, k$ remain fixed).
\begin{proof}   Corollary \ref{cor:Uniqueness} and Lemma \ref{lem:Ccount} imply \[|F+C(m,n;k)|=|F||C(m,n;k)|=|F|\Bigl(\frac{1}{|G_{m}|}+o(1)\Bigr)|G_{n}|. \qedhere\]
\end{proof}

Considering the restriction of an element of $C(m,n;k)$ to an interval $\tau\in \Omega_{r}$, an argument similar to the proof of Lemma \ref{lem:Vrestriction} yields the following.
\begin{observation}\label{obs:Crestriction}
  If $r\leq m<n$, $k\in \mathbb N$ and $g\in C(m,n;k)$, then for all $\tau\in \Omega_{r}$, $g_{\tau}\in C(m-r,n-r;k)$.
\end{observation}

If $S_{1},\dots, S_{l}\subseteq \mathbb F_{2}^{\omega}$ is a sequence of sets, write $\sum_{i=1}^{l} S_{i}$ for the set $\{s_{1}+\cdots +s_{l}:s_{i}\in S_{i} \text{ for all }i \leq l\}$.

\begin{lemma}\label{lem:Cshift}
  Let $0=n_{0}<n_{1}<n_{2}<\cdots<n_{l}$ be an increasing sequence of integers and  $k_{i}, k_{i}'>0$  for each $i$. If $A:=\sum_{i=1}^{l} C(n_{i-1},n_{i};k_{i})$, $A':=\sum_{i=1}^{l} C(n_{i-1},n_{i};k_{i}')$, then $(A+\mathbf{1})\cap A'=\varnothing$.
\end{lemma}

\begin{proof}  We proceed by induction on $l$.  For $l=1$, the conclusion follows from Lemma \ref{lem:BaseDisjoint}.  Now suppose $l>1$, and the conclusion holds for every increasing sequence $0 = n_{0}'< n_{1}' < \cdots <n_{l-1}'$.
Assume, to get a contradiction, that the intersection is nonempty.  So there are $c_{i}\in C(n_{i-1},n_{i};k_{i})$, $c_{i}'\in C(n_{i-1},n_{i};k_{i}')$ such that
\begin{equation}\label{eqn:CC'}
  c_{1}+\cdots + c_{l}+\mathbf 1= c_{1}'+ \cdots + c_{l}'.
\end{equation}
Let $h=c_{1}$, $h'=c_{1}'$, $f=c_{2}+\cdots + c_{j}$, $f'=c_{2}'+\cdots + c_{j}'$.   Then $h \in C(0, n_{1};k_{1})$, $h'\in C(0,n_{1};k_{1}')$. Since $\|h\|_{0}+\|h'\|_{0}>1$, there exists $t\in I$ such that $h(t)=h'(t)=0$.  Since both $h$, $h'$ are constant on elements of $\Omega_{n_{1}}$, there is a  $\tau\in \Omega_{n_{1}}$ such that $h_{\tau}=h_{\tau}'=\mathbf 0$.  For such $\tau$, the definition of $f$, $f'$ and Observation \ref{obs:Crestriction} yield
\begin{align*}
f_{\tau} \in \sum_{i=2}^{l} C(n_{i-1}-n_{1}, n_{i}-n_{1}, k_{i}), &&  f_{\tau}' \in \sum_{i=2}^{l} C(n_{i-1}-n_{1}, n_{i}-n_{1}, k_{i}'),
\end{align*}
while Equation (\ref{eqn:CC'}) implies $f_{\tau}+\mathbf{1}=f_{\tau}'$.  The induction hypothesis says that $f_{\tau}+\mathbf{1}\neq f_{\tau}'$, so we have the desired contradiction.  \end{proof}

\begin{lemma}\label{lem:CplusU}
Let $m\leq n, k'\leq  k \in \mathbb N\cup \{0\}$.  Then
\begin{enumerate}
  \item[(i)] $C(m,n;k)+U(n,k')\subseteq C(m,n;k-k')$,
\item[(ii)]  $C(m,n;k)\subseteq C(m,n;k')$.
 \end{enumerate}

\end{lemma}

\begin{proof}
(i)  Lemma \ref{lem:Vrestriction} implies that for all $u\in U(n,k')$ and all $\tau \in \Omega_{m}$, we have $u_{\tau}\in U(n-m,k')$.  If $h\in C(m,n;k)$, then $h_{\tau}\in C(0,n-m,k)$. Lemma \ref{lem:uShift} implies $h_{\tau}+u_{\tau} \in C(0,n-m,k-k')$ for all $\tau\in \Omega_{m}$, and the conclusion now follows from the definition of $C(m,n;k-k')$.

Part (i) follows from part (ii) of Lemma \ref{lem:uShift} and the definition of $C(m,n;k)$.
\end{proof}

\subsection{Proofs.}
In this section we prove Theorems \ref{thm:Main} and \ref{thm:PWsynd}.

Fix $\varepsilon>0$.    Apply Lemma \ref{lem:BaseCardinality} to choose $n_{1}$ sufficiently large that $\frac{|C(0,n_{1};3)|}{|G_{n_{1}}|}>\frac{1}{2}-\frac{\varepsilon}{8}$.  Apply Lemma \ref{lem:Cardinality} to choose $n_{2}$ sufficiently large that \[
\frac{|C(0,n_{1};3)+C(n_{1},n_{2};6)|}{|G_{n_{2}}|}>\Bigl(1-\frac{\varepsilon}{4}\Bigr)\frac{|C(0,n_{1};3)|}{|G_{n_{1}}|}>\frac{1}{2}-\frac{\varepsilon}{2}.
\] Continuing in this way,  apply Lemma \ref{lem:Cardinality} to choose $n_{3}<n_{4}<\cdots$ sufficiently large that
\begin{equation}\label{eqn:Clarge}
 \frac{|C(0,n_{1};3)+\cdots + C(n_{l-1},n_{l};3l)|}{|G_{n_{l}}|} > \frac{1}{2} - \varepsilon
\end{equation}
for every $l$.  Let
\begin{align*}
  A_{l}&:=C(0,n_{1};3)+C(n_{1},n_{2};6)+\cdots + C(n_{l-1},n_{l};3l),\\
A_{l}'&:=C(0,n_{1};2)+C(n_{1},n_{2};4)+\cdots + C(n_{l-1},n_{l};2l).
\end{align*}
Note that
\begin{eqnarray}\label{eqn:Monotone}
   A_{l}\subseteq A_{l+1},\ A_{l}'\subseteq A_{l+1}' && \text{ for all } l,
\end{eqnarray} since $\mathbf 0 \in C(n_{l},n_{l+1};k)$ for each $k$.

Let $A:=\bigcup_{l=1}^{\infty}A_{l}$, $A':=\bigcup_{l=1}^{\infty} A_{l}'$.
Let $T:=\bigcup_{j=1}^{\infty} U(n_{j},j)$, and let $S:=T+\mathbf 1$. We claim that \begin{equation}\label{eqn:AA'}
  (A+T) \subseteq A'.
\end{equation}
If we show $A_{l}+U(n_{j},j)\subseteq A_{l}'$ whenever $j\leq l$ then the containment (\ref{eqn:Monotone}) yields (\ref{eqn:AA'}).  Part (i) of Lemma \ref{lem:CplusU} implies $U(n_{j},j)+C(n_{j-1},n_{j};3j)\subseteq C(n_{j-1},n_{j};2j)$, and Part (ii) of Lemma \ref{lem:CplusU} yields the desired containment.

To prove Theorem \ref{thm:Main}, it suffices to prove the following.
\begin{eqnarray}\label{eqn:Adensity}
  &d^{*}(A)\geq \tfrac{1}{2}-\varepsilon,\\
  &\label{eqn:Disjoint} (A+S)\cap A=\varnothing,\\
\label{eqn:gS} &g+S \text{ is chromatically intersective for all } g\in G.
\end{eqnarray}
Inequality (\ref{eqn:Clarge}) implies $\frac{|A_{l}|}{|G_{n_{l}}|}\geq \frac{1}{2}-\varepsilon$ for each $l$.  Observation \ref{obs:UBD}   then implies Inequality (\ref{eqn:Adensity}).

To prove (\ref{eqn:Disjoint}),  note that $A+S=A+T+\mathbf 1 \subseteq A'+\mathbf 1$, by (\ref{eqn:AA'}). Lemma \ref{lem:Cshift} implies that $(A_{l}'+\mathbf{1})\cap A_{l}=\varnothing$ for every $l$, hence $(A'+\mathbf{1})\cap A=\varnothing$.   Then $(A+S)\cap A\subseteq (A'+\mathbf 1)\cap A=\varnothing$.

Lemma \ref{lemmaChromatic} implies (\ref{eqn:gS}).  Consequently, $(g+B-B)\cap S\neq \varnothing$ whenever $g\in G$ and $B\subseteq G$ is piecewise syndetic.  Equation (\ref{eqn:Disjoint}) implies $(A-A)\cap S=\varnothing$, hence $A-A$ does not contain such a $g+B-B$. We have proved Theorem \ref{thm:Main}.

To prove Theorem \ref{thm:PWsynd}, let $E:=A\cup (A+\mathbf 1)$.  Equation (\ref{eqn:Disjoint}) implies  $A$ and $A+\mathbf{1}$ are disjoint, so Inequality (\ref{eqn:Adensity}) and Lemma \ref{lemAdditivity} imply $d^{*}(E)>1-2\varepsilon$.   It therefore suffices to show that $E+S$ is not piecewise syndetic. Setting $D:=A+T$, (\ref{eqn:AA'}) implies $D-D\subseteq A'-A'$.  We may repeat the proof that $(A-A)\cap S=\varnothing$ to show that $(A'-A')\cap S=\varnothing$.  Since $(B-B)\cap S\neq \varnothing$ whenever $B$ is piecewise syndetic, it follows that $A'$ is not piecewise syndetic, meaning $A+T$ is not piecewise syndetic.   The partition regularity of piecewise syndeticity (see \cite{HindmanStrauss}) and the identity $E+S=(A+T+\mathbf 1)\cup (A+T)$ imply that $E+S$ is not piecewise syndetic. This completes the proof of Theorem \ref{thm:PWsynd}.

\section{Proof of Theorems \ref{thm:Oddp} and \ref{thm:PWBohr}}\label{sec:Oddp}
  When $p$ is odd our construction is similar to the proof of Theorem \ref{thm:Main}, but the construction of dense sets is more intricate.

\subsection{Presentation of \texorpdfstring{$\mathbb F_{p}^{\omega}$}{FPw}} Fix a prime number $p$.  Write the elements of $\mathbb F_{p}$ as $0, 1, \dots, p-1$.  We continue to write $I$ for the half-open unit interval $[0,1)$, $\mu$ for Lebesgue measure on $I$, and $\Omega_{n}$ for the collection of intervals $\{[\frac{j}{2^{n}}, \frac{j+1}{2^{n}}):0\leq j \leq 2^{n}-1\}$.  Let $G_{n}$ be the group of functions $f:I\to \mathbb F_{p}$ which are constant on every interval $\tau\in \Omega_{n}$, with the group operation of pointwise addition.  Then $G_{n}$ is isomorphic to $\mathbb F_{p}^{2^{n}}$, and has cardinality $p^{2^{n}}$. Let $G=\bigcup_{n=1}^{\infty} G_{n}$.  Then $G$ is a countably infinite vector space over $\mathbb F_{p}$, so $G$ is isomorphic to $\mathbb F_{p}^{\omega}$.

Let $\mathbf 0$  denote the identity element of $G$, $\mathbf 1$ denote the constant function $\mathbf 1(t)=1\in \mathbb F_{p}$ for all $t\in I$, and for $y\in \mathbb F_{p}$, let $y\mathbf 1$ denote the constant function so that $y\mathbf 1(t)=y$ for all $t\in I$.

For $g\in G$ and $x\in \mathbb F_{p}$, let $\|g\|_{x}:=\mu(g^{-1}(x))$. The following identity is crucial.
\begin{equation}\label{eqn:WeightShift}
  \|g+y\mathbf 1\|_{x} = \|g\|_{x-y} \text{ for all } g\in G, x,y\in \mathbb F_{p}.
\end{equation}
The identity follows from the equation $(g+y\mathbf{1})^{-1}(x)=g^{-1}(x-y)$.

For $k,n \in \mathbb N$, let $V(n,k):=\{g\in G_{n}:\|x\|_{1}\geq1-k2^{-n}\}$.

\begin{remark}
  We do not assume $p$ is odd here, and in fact when $p=2$ the construction of the sets $A$ and $S$ in this section will repeat the construction of $A$ and $S$ in \S\S \ref{sec:Groups} and \ref{section2Proof}.
\end{remark}

\begin{lemma}\label{lem:BohrIntersective}
  Let $G$ be a countable abelian group, $g\in G$, $S\subset G$ dense in the Bohr topology, and $B\subset G$ a piecewise Bohr set.  Then $g+(B-B)\cap S\neq \varnothing$.
\end{lemma}
\begin{proof}
  There is a Bohr set $B'$ such that $B$ contains a translate of every finite subset of $B'$.  Hence $B'-B'\subset B-B$.  Since $B'$ is a Bohr set, $g+(B'-B')$ is a Bohr neighborhood, so it has nonempty intersection  with $S$.
\end{proof}

\subsection{Bohr topology of \texorpdfstring{$\mathbb F_{p}^{\omega}$}{FPw} and Hamming balls}\label{sec:pHammingBalls}  Since every nonzero element of $\mathbb F_{p}^{\omega}$ has order $p$, every Bohr set is a finite intersection of kernels of homomorphisms $\rho:\mathbb F_{p}^{\omega}\to Z_{p}$, where $Z_{p}\subseteq \mathbb C$ is the group of $p$th roots of unity.  It follows that the Bohr topology on $\mathbb F_{p}^{\omega}$ is the topology whose open sets are unions of cosets of finite index subgroups.  This leads to the following observation.

\begin{observation}\label{obs:BohrFinite}
  A set $S\subseteq \mathbb F_{p}^{\omega}$ is dense in the Bohr topology if and only if $\rho(S)=\rho(G)$ for every surjective homomorphism $\rho$ from $G$ to a finite group.
\end{observation}
\begin{observation}\label{obs:Generate}
  If $E\subseteq \mathbb F_{p}^{k}$ generates $\mathbb F_{p}^{k}$ as an abelian group, then every element of $\mathbb F_{p}^{k}$ is a sum of at most $pk$ elements of $E$, since some subset of $E$ must form a basis of $\mathbb F_{p}^{k}$ as a vector space over $\mathbb F_{p}$.
\end{observation}

 In order to prove that unions of the $V(n,k)$ are dense in the Bohr topology, we first consider Hamming balls in $\mathbb F_{p}^{n}$.  For $x\in \mathbb F_{p}^{n}$, let $\supp(x):=\{i\leq n: x_{i}\neq 0\}$, and let $|x|$ be the cardinality of $\supp(x)$.  Let $H_{n,k}:=\{x\in \mathbb F_{p}^{n}:|x|\leq k\}$. We say that $S\subseteq \mathbb F_{p}^{n}$ is \emph{$k$-Bohr dense} if for every surjective homomorphism $\rho: \mathbb F_{p}^{n}\to \mathbb F_{p}^{k}$, $\rho(S)=\mathbb F_{p}^{k}$.

\begin{lemma}\label{lem:kBohrDense}
$H_{n,pk}$ is $k$-Bohr dense in $\mathbb F_{p}^{n}$.
\end{lemma}

\begin{proof}
  Let $\rho: \mathbb F_{p}^{n}\to \mathbb F_{p}^{k}$ be a surjective homomorphism.    Now $H_{n,1}$ generates $\mathbb F_{p}^{n}$, so $\rho(H_{n,1})$ generates $\mathbb F_{p}^{k}$, and Observation \ref{obs:Generate} guarantees that every element of $\mathbb F_{p}^{k}$ is a sum of at most $pk$ elements of $\rho(H_{n,1})$.  In other words, $\mathbb F_{p}^{k}=\sum_{j=1}^{pk} \rho(H_{n,1})$. Since $\rho$ is a homomorphism, the last set is equal to $\rho(H_{n,pk})$.
\end{proof}

\begin{lemma}\label{lem:pBorhDense}
  If $n_{j}$, $k_{j}\to \infty$, then $T:=\bigcup_{j=1}^{\infty} U(n_{j},k_{j})$ is dense in the Bohr topology of $G$.
\end{lemma}
\begin{proof}
By Observation \ref{obs:BohrFinite} it suffices to show that $\rho(T)=K$ for every surjective homomorphism $\rho$ from $G$ to a finite group $K$.  Let $\rho$ be such a homomorphism. Then every nonzero element of $K$ has order $p$, so $K$ is isomorphic to $\mathbb F_{p}^{m}$ for some $m\in \mathbb N\cup \{0\}$.  Choose $j$ so that $\rho(G_{n_{j}})=K$ and $k_{j}>pm$.  The natural isomorphism of $G_{n_{j}}$ with $\mathbb F_{p}^{2^{n_{j}}}$ identifies $U(n_{j},k_{j})$ with $H_{2^{n_{j}},k_{j}}$, so Lemma \ref{lem:kBohrDense} implies $\rho(U(n_{j},k_{j}))=K.$
\end{proof}

\subsection{Bias patterns}

Let $\mathcal S:=\{S\subseteq \mathbb F_{p}: S\notin \{\varnothing, \mathbb F_{p}\}\}$.  Then $\mathbb F_{p}$ acts on $\mathcal S$ by addition: $S+x:=\{s+x: s\in S\}$.  The orbit of every $S$ has cardinality $p$, since $p$ is prime and $S+1\neq S$ for all $S\in \mathcal S$.  Fix, once and for all, a partition of $\mathcal S$ into sets $\mathcal S_{x}$, $x\in \mathbb F_{p}$, such that $\mathcal S_{x+y}=\{S+y:S\in \mathcal S_{x}\}$ for all $x, y\in \mathbb F_{p}$.  We also insist that $\{0\}\in \mathcal S_{0}$.  For example, with $p=3$, we can take
\begin{align*}
  \mathcal S_{0} =\{\{0\},\{0,1\}\},\ \mathcal S_{1} =\{\{1\},\{1,2\}\},\ \mathcal S_{2} =\{\{2\},\{2,0\}\}.
\end{align*}

For $S\in \mathcal S$ and $n,k\in \mathbb N$ let $\Bias_{n}(S,k)$ be the set of $g\in G_{n}$ such that
\begin{align}
\label{eqn:B1}  \|g\|_{x}&\geq  \tfrac{1}{p}+k2^{-n} && \text{for all } x\in S\\
\label{eqn:B2}  \|g\|_{x} &\leq  \tfrac{1}{p}-k2^{-n} && \text{for all } x\in \mathbb F_{p}\setminus S.
\end{align}
\begin{lemma}\label{lem:BiasBase} Let $n,k \in \mathbb N$, $y\in \mathbb F_{p}$.
\begin{enumerate}
  \item[(i)]  $\Bias_{n}(S,k)+y\mathbf 1 = \Bias_{n}(S+y,k)$.
\smallskip
 \item[(ii)] If $k'>0$, $S\neq S'\in \mathcal S$, then $\Bias_{n}(S,k)\cap \Bias_{n}( S', k')=\varnothing$. \end{enumerate}
\end{lemma}
\begin{proof} (i) We will show that $\Bias_{n}(S,k)+y\mathbf 1 \subseteq \Bias_{n}(S+y,k)$; the reverse containment follows by symmetry.  Suppose $g\in \Bias_{n}(S,k)$ and $y\in \mathbb F_{p}$. If $x\in S+y$, then $x-y\in S$, and Identity (\ref{eqn:WeightShift}) implies $\|g+y\mathbf 1\|_{x}=\|g\|_{x-y}\geq \frac{1}{p}+k2^{-n}$.  Similarly, if $x\in \mathbb F_{p}\setminus (S+y)$, then $\|g+y\mathbf 1\|_{x}\leq \frac{1}{p}-k2^{-n}$.  Thus $g+y\mathbf 1\in \Bias_{n}(S+y,k)$.

To prove Part (ii), let $g$ be in the intersection.  If $S\neq S'$, there is an $x$ in the symmetric difference $S\triangle S'$. For such $x$, Inequalities (\ref{eqn:B1}) and (\ref{eqn:B2}) must both be satisfied, which is impossible.
\end{proof}

\begin{definition}\label{def:pC0}
  For $n,k\in \mathbb N$, let $C(0,n;k):=\bigcup_{S\in \mathcal S_{0}} \Bias_{n}(S,k)$.
\end{definition}

Note that our insistence that $\{0\}\in \mathcal S_{0}$ implies $\mathbf 0\in C(0,n;k)$.

\begin{lemma}\label{lem:CpDisjoint}
  If $x\in \mathbb F_{p}\setminus \{0\}$ and $k,k'>0$, then $(C(0,n;k)+x\mathbf 1)\cap C(0,n;k')=\varnothing$.
\end{lemma}
\begin{proof}  It suffices to show that $(\Bias_{n}(S,k)+x\mathbf 1)\cap \Bias_{n}(S',k')=\varnothing$  for every $S,S'\in \mathcal S_{0}$.  Fix such $S$, $S'$.  Part (i) of Lemma \ref{lem:BiasBase} implies $\Bias_{n}(S,k)+x\mathbf 1=\Bias_{n}(S+x,k)$.  Now $S+x\in \mathcal S_{x}$, which is disjoint from $\mathcal S_{0}$ by definition.  Hence $S+x\neq S'$, and Part (ii) of Lemma \ref{lem:BiasBase} provides the desired disjointness.
\end{proof}

\begin{lemma}\label{lem:BasepCard}
  Fix $k\in \mathbb N$.  Then $\lim_{n\to \infty} \frac{|C(0,n;k)|}{|G_{n}|}=\frac{1}{p}$.
\end{lemma}
\begin{proof}
   For $x\in \mathbb F_{p}$, let $E_{n,x}:=C(0,n;k)+x\mathbf 1$, and let $Z_{n}:=G_{n}\setminus \bigcup_{x\in \mathbb F_{p}} E_{n,x}$. By Lemma \ref{lem:CpDisjoint}, the sets $E_{n,x}$, $x\in \mathbb F_{p}$ are mutually disjoint, and $|E_{n,x}|=|E_{n,y}|$ for all $x,y\in \mathbb F_{p}$.  It therefore suffices to show that $|Z_{n}|=o(|G_{n}|)$. Note that $Z_{n}$ is the set of $g\in G_{n}$ such that $\frac{1}{p}-k2^{-n}<\|g\|_{x}<\frac{1}{p}+k2^{-n}$ for every $x\in \mathbb F_{p}$.

    Let $M_{n,k}$ be the largest value of $\binom{|\Omega_{n}|}{t}$ where  $\frac{1}{p}|\Omega_{n}| - k \leq  t \leq \frac{1}{p}|\Omega_{n}| + k$.  For such $t$  and $y\in \mathbb F_{p}$, the number of $g\in G_{n}$ satisfying $\|g\|_{y}= t2^{-n}$ is at most
\[
   (p-1)^{|\Omega_{n}|-\lfloor\frac{1}{p}|\Omega_{n}|-k\rfloor}M_{n,k}.
\]
Summing over all possible values of $t$ and $y$, we get
\[
|Z_{n}| \leq p(2k+1)(p-1)^{|\Omega_{n}|-\lfloor\frac{1}{p}|\Omega_{n}|-k\rfloor}M_{n,k}.
\]Estimating the binomial coefficients in the definition of $M_{n,k}$ with Stirling's formula, we find $|Z_{n}| = o(p^{|\Omega_{n}|}) = o(|G_{n}|)$. \end{proof}

\begin{definition}\label{def:PC}
  For $m\in \mathbb N\cup \{0\}$, $n>m$ and $k\in \mathbb N$, let $C(m,n;k)$ be the set of $g\in G_{n}$ such that $g_{\tau}\in C(0,n-m;k)$ for every $\tau\in \Omega_{m}$.
\end{definition}

An argument similar to the proof of Lemma \ref{lem:Ccount} shows that
\begin{equation}\label{eqn:pCcount}
  |C(m,n;k)|=\Bigl(\frac{1}{|G_{m}|}+o(1)\Bigr)|G_{n}|
\end{equation}
as $n\to \infty$ and $m,k$ remain fixed.

\begin{lemma}\label{lem:pShift}
  With $m,n,k$ as in the above definition, if $g, g' \in G_{m}$ and $h,h'\in C(m,n;k)$, then $g+h=g'+h'$ if and only if $g=g'$ and $h=h'$.
\end{lemma}

\begin{proof}
  If $g+h=g'+h'$, then $g-g'=h-h'$.  Since $g-g'\in G_{m}$, we have that $(g-g')_{\tau}=x_{\tau}\mathbf{1}$ for every $\tau\in \Omega_{m}$, where $x_{\tau}\in \mathbb F_{p}$.  Then $(h-h')_{\tau} =x_{\tau}\mathbf{1}$ for every $\tau\in \Omega_{m}$, meaning
  \begin{equation}\label{eqn:CellShift}
  h_{\tau}\in (C(0,n-m;k)+x_{\tau}\mathbf{1})\cap C(0,n-m;k)
  \end{equation}
  for all $\tau \in \Omega_{m}$.  Lemma \ref{lem:pShift} and the inclusion (\ref{eqn:CellShift}) imply $x_{\tau}=0$, so $g-g'=\mathbf 0$ and $h-h'=\mathbf 0$.
\end{proof}

\begin{corollary}\label{cor:Psumset}
  With $n$, $m$, and $k$ as in Definition \ref{def:PC} and $F\subseteq G_{m}$, $|F+C(n,m;k)|=|F||C(n,m;k)|$.
\end{corollary}

The remainder of the proof of Theorems \ref{thm:Oddp} and \ref{thm:PWBohr} follows steps analogous to Corollary \ref{cor:Uniqueness} and the subsequent remainder of Section \ref{section2Proof}.  The proof of the following lemma is exactly analogous to the proof of Lemma \ref{lem:Cardinality}, using Lemma \ref{lem:pCard} in place of Lemma \ref{lem:Ccount} and Corollary \ref{cor:Psumset} in place of Corollary \ref{cor:Uniqueness}.

\begin{lemma}\label{lem:pCard}
If $F\subseteq G_{m}$ then $\lim_{n\to \infty} \frac{|F+C(m,n;k)|}{|G_{n}|}=\frac{|F|}{|G_{m}|}$.
\end{lemma}

\begin{lemma}\label{lem:pCshift}
    Let $0=n_{0}<n_{1}<n_{2}<\cdots<n_{l}$ be an increasing sequence of integers and  $k_{i}, k_{i}'>0$  for each $i$. If $A:=\sum_{i=1}^{l} C(n_{i-1},n_{i};k_{i})$, $A':=\sum_{i=1}^{l} C(n_{i-1},n_{i};k_{i}')$ and $x\in \mathbb F_{p}\setminus \{0\}$, then $(A+x\mathbf{1})\cap A'=\varnothing$.
\end{lemma}

\begin{proof}
  Induction on $l$.  The base case $l=1$ is Lemma \ref{lem:CpDisjoint}.  Assume the lemma holds for the sequence $0=n_{0}<n_{1}<\cdots<n_{l-1}$ and $k_{i}$, $k_{i}'>0$.  Assume, to get a contradiction, that $(A+x\mathbf 1)\cap A \neq\varnothing$.  Then there are $c_{i}\in C(n_{i-1},c_{i};k_{i})$, $c_{i}'\in C(n_{i-1},n_{i};k_{i}')$ such that
  \begin{equation}\label{eqn:Crearrange}
   x\mathbf{1}+c_{1}+\cdots+c_{l} = c_{1}'+\cdots + c_{l}'
  \end{equation}
  The left-hand side of Equation (\ref{eqn:Crearrange}) is an sum of the form $g+h$, where $g=x\mathbf{1}+c_{1}+\cdots+c_{l-1} \in G_{n_{l-1}}$ and $h=c_{l}\in C(n_{l-1},n_{l};k)$, and the right hand side has a similar form, as $g'+h'$, where $g'=c_{1}'+\cdots+c_{l-1}'\in G_{n_{l-1}}$ and $h'=c_{l}' \in C(n_{l-1},n_{l};k')$.  Both $c_{l}$ and $c_{l}'$ are in $C(n_{l-1},n_{l};k'')$, where $k''=\min(k,k')$.  Lemma \ref{lem:pShift} then implies $g=g'$, so $x\mathbf 1+c_{1}+\cdots + c_{l-1}= c_{1}'+\cdots +c_{l-1}'$.  The last equation is impossible, by the induction hypothesis, and this is the desired contradiction. \end{proof}

\subsection{Proof of Theorems \ref{thm:Oddp} and \ref{thm:PWBohr}}  Fix $\varepsilon>0$.  Apply Lemma \ref{lem:BasepCard} to choose $n_{1}$ large enough that $|C(0,n_{1},3)|>\frac{1-\varepsilon}{p}|G_{n_{1}}|$.   Apply Lemma \ref{lem:pCard} to choose a sequence of integers $n_{2}<n_{3}<\cdots$ so that
\begin{equation}\label{eqn:pSumCard}
  |C(0,n_{1},3)+\cdots + C(n_{l-1},n_{l},3l)| > \frac{1-\varepsilon}{p} |G_{n_{l}}|
\end{equation}
for each $l$.  Set $B_{l}:=\sum_{j=1}^{l} C(n_{j-1}, n_{j};3j)$.  Let $R:=\{0,2,\dots, p-3\}\subseteq \mathbb F_{p}$, so that $|R|=\frac{p-1}{2}$.  Let
\[
A_{l}:=\bigcup_{x\in R} B_{l}+x\mathbf 1.
\]
Finally, let $A:=\bigcup_{l=1}^{\infty} A_{l}$.  The sets $B_{l}+x\mathbf{1}$, $x\in \mathbb F_{p}$ are mutually disjoint, by Lemma \ref{lem:pCshift}, so $d^{*}(A)\geq |R| \frac{1-\varepsilon}{p}=\frac{p-1}{2p}(1-\varepsilon)$.  Thus  $d^{*}(A)\geq \frac{1}{2}-\frac{1}{2p}-\varepsilon$.

Let $S:=\bigcup_{j=1}^{\infty} V(n_{j},j)$.  Lemma \ref{lem:pBorhDense} implies $S$ is dense in the Bohr topology of $G$.

Lemma \ref{lem:pCshift} implies that $(A+S)\cap A=\varnothing$, via an argument similar to the proof of (\ref{eqn:Disjoint}), so $(A-A)\cap S=\varnothing$, proving Theorem \ref{thm:Oddp}.   Likewise we have that $A':=A+S$ satisfies $(A'-A')\cap S=\varnothing$.  Thus $A'$ cannot be piecewise Bohr, as Lemma \ref{lem:BohrIntersective} guarantees $(B-B)\cap S\neq \varnothing$ for every piecewise Bohr set $B$.  This completes the proof of Theorem \ref{thm:PWBohr}.

\subsection{Questions}
If $G$ is an abelian group and $S\subseteq G$, let $\chi(S)$ be the greatest $m$ such that $S$ is $(m-1)$-chromatically intersective, or $\chi(S)=\infty$ if the set of such $m$ is unbounded.

\begin{question}\label{question:ChromaticNumber}
  Fix an odd prime $p$, and let $H_{n,k}$ be defined as in \S\ref{sec:pHammingBalls}.  Does $\chi(g+H_{n,k})\to \infty$ as $k\to\infty$, independently of $n$ and $g\in \mathbb F_{p}^{n}$?
\end{question}
If the answer to Question \ref{question:ChromaticNumber} is ``yes," then the conclusions Theorems \ref{thm:Oddp} and \ref{thm:PWBohr} can be improved to match the conclusions of Theorems \ref{thm:Main} and \ref{thm:PWsynd}.

\begin{question}\label{question:Heredetary}
  Let $G$ be a countable abelian group. Which, if any, of the following implications hold?

\begin{enumerate}
  \item[(I)] If every translate of $S$ is a set of topological recurrence, then there is a set $S'\subseteq S$ such that every translate of $S'$ is a set of topological recurrence while for all $g\in G$, $(S'-g)\setminus \{0\}$ is not a set of measurable recurrence.

\item[(II)] If $S\subseteq G$ is a set of topological recurrence, there a set $S'\subseteq S$ such that $S'$ is a set of topological recurrence but not a set of measurable recurrence.

\item[(III)]  If $S$ is dense in the Bohr topology of $G$, there is a set $S'\subseteq S$ such that $S'$ is dense in the Bohr topology of $G$ while for all $g\in G$, $(S'-g)\setminus \{0\}$ is not a set of measurable recurrence. \end{enumerate}
\end{question}

See \cite{GrSingle} for a list of problems related to our results.

\begin{remark}\label{rem:Presentation}  We have two reasons for using the presentation of $\mathbb F_{2}^{\omega}$ introduced in \S\ref{sec:Groups}.  The first is that we do not know how to prove Theorem \ref{thm:Main} using the standard presentation.  To elaborate: write an element of $\mathbb F_{2}^{\omega}$ as $x=(x_{1},x_{2},\dots,)$, and let $|x|=|\{i:x_{i}=1\}|$.  Let $F_{n}:=\{x\in \mathbb F_{2}^{\omega}:x_{i}=0 \text{ for all } i>n \}$. The natural approach to proving Theorem \ref{thm:Main} would be: let $(n_{j})_{j\in \mathbb N}$ be a rapidly increasing sequence of integers, let $S_{j} \subseteq \mathbb F_{2}^{\omega}$, $S_{j}:=\{x\in F_{n_{j}}: |x|>n_{j}-j\}$.  Then Lemma \ref{lem:ChromaticBase} implies every translate of $S:=\bigcup_{j=1}^{\infty} S_{j}$ is chromatically intersective.  One would then construct sets $A_{l} \subseteq F_{n_{l}}$ having $|A_{l}|\approx \frac{1}{2}|F_{n_{l}}|$ and $(A_{l}-A_{l})\cap S=\varnothing$ for each $l$. We do not see how to construct such $A_{l}$, so we resorted to the presentation of $\mathbb F_{2}^{\omega}$ given in \S\ref{sec:Groups}.

Our second reason for using our preferred presentation is that it arises naturally in certain approaches to answer Question \ref{question:Heredetary} for $G=\mathbb F_{p}^{\omega}$.
\end{remark}

  \bibliographystyle{amsplain}
\frenchspacing
\bibliography{Bohr_recurrence_revision2}

\end{document}